\documentclass[12pt]{amsart}
\usepackage{amsfonts,amsmath,amssymb,eucal,tipa}
\usepackage[all]{xy}
\usepackage{hyperref}
\usepackage[utf8]{inputenc}

\begin{document}

\parindent=0pt
\parskip=6pt

\newcommand{\SP}{{\rm SP}}
\newcommand{\B}{{\mathcal B}}
\newcommand{\C}{{\mathbb C}}
\newcommand{\E}{{\mathbb E}}
\newcommand{\ka}{{\kappa}}
\newcommand{\G}{{\mathbb {G}_m}}
\newcommand{\hu}{{\mathfrak h}}
\newcommand{\ve}{{\mathfrak e}}
\newcommand{\Q}{{\mathbb Q}}
\newcommand{\CP}{{\mathbb{C}P}}
\newcommand{\MU}{{\rm MU}}
\newcommand{\PGl}{{\rm PGl}}
\newcommand{\R}{{\mathbb R}}
\newcommand{\Sf}{{\sf S}}
\newcommand{\SU}{{\rm SU}}
\newcommand{\T}{{\mathbb T}}
\newcommand{\te}{{\hbar \nu/kT}}
\newcommand{\N}{{\mathbb N}}
\newcommand{\st}{{\sf st}}
\newcommand{\Z}{{\mathbb Z}}

\newcommand{\tkap}{{\circ \kappa^{-1} \circ}}
\newcommand{\tsf}{{\circ \sf{st}^{-1} \circ}}

\title{Notes toward a Newtonian thermodynamics}

\author[Jack Morava]{Jack Morava}

\address{Department of Mathematics, The Johns Hopkins University,
Baltimore, Maryland}

\begin{abstract}{We interpret the moment generating function 
\[
{\mathbb E}(\exp_\G(tX)) := \exp_F(t) \in \R[[x]]
\]
of a random variable $X$ as the exponential of an associated one-dimensional formal group law $F$ defined over $\R$.}
\end{abstract}

\maketitle \bigskip

{\bf Introduction} In $p$-adic Fourier theory (\cite{11}(lemma 3.4, diagram at foot of p 460) the Schneider-Teitelbaum modulus\begin{footnote}{the natural logarithm if $X$ is Poisson; the Bernoulli distribution yields Todd's arithmetic genus}\end{footnote}
\[
\st_F = \exp_\G \circ \log_F
\]
of $F$ and the cumulant generating function 
\[
\kappa_F = \log_\G \circ \exp_F 
\]
of $X$ (roughly, its associated Gibbs free energy \cite{12}) are, up to intertwining coordinate transformations
\[
\kappa_F = \log_\G \tsf \exp_\G, \;  \st_F = \exp_\G \tkap \log_\G \;,
\]
mutually inverse as functions. Quillen's interpretation of the complex cobordism ring $\MU_*$ as classifying formal groups expresses $\st$, and therefore $\kappa$, in terms of universal symmetric symmetric functions with rational coefficients. The symplectic cobordism ring $\B_*$  of VL Ginzburg  \cite{3,7} provides a thermodynamic extension of this interpretation in which 
\[
[\CP_1,\omega] \sim 1/kT \cdot [\CP_1]  \dots 
\]
\bigskip

{\bf \S 1.1} Recall, re symmetric functions \cite{5}(I \S 2), that
\[
E(t) = \prod (1 + x_i t) , \; H(t) = E(-t)^{-1} 
\]
\[
P(t) = H'(t)/H(t) = \sum p_r t^r \; \Rightarrow \; H(t) = \exp(p_1 \cdot \sum \wp_n \frac{t^n}{n}) \;,
\]
with slightly eccentric notation $p_n = p_1 \cdot \wp_n$ for the power sums $p_n = \sum x^n_i$.
These are elements of a periodically-graded algebra $\Sf [e_1^{-1}]$ of symmetric functions (over $\Z$), with one generator in each even degree.

[This algebra has at least a Klein four-group $\{e,h;m,f\}$ of useful bases \cite{9}(\S 1.2), corresponding to the operations on virtual vector bundles (and their characteristic classes) defined by both inverse $V \mapsto - V$ and complex conjugation $V \mapsto \overline {V}$; the semicolon indicates duality with respect to the Hall (positive definite integral) inner product \cite{5}(I \S 4.9). Basis change can be expressed in terms of Kostka matrices \cite{5}(I \S 6).] \bigskip

{\bf \S 1.2} Recall the level one projective representation 
\[
\Omega \SU(2) \to \PGl_\C({\mathbb H}) \; ;
\]
then

* ) $\; \CP_\infty \simeq B\T \simeq H(\Z,2)$ so $\pi_2(\CP_\infty) \cong \Z$
 
**)  $\; \Omega \SU(2) \simeq (\Omega \Sigma) \;  S^2 \simeq \vee_{n \geq 0} S^{2n}$ stably

***) $b = b_{(1)} = [S^2 = \CP_1 \subset \CP_\infty]$ generates $\pi_2$. Let 
\[
b_{(n)} = [\CP_n \subset \CP_\infty] \in \MU_{2n}\CP_\infty
\]
as in \cite{11}.

It follows from **)  that $H_* \Omega \SU(2) \cong \Z [b]$ as Pontrjagin rings, while $H^*\CP_\infty \cong \Z[c]$ as Hopf algebras. Let us write $H_*\CP_\infty \cong \Z[\gamma^*b]$ with divided power operations $\gamma^k x = x^k/k!$ and $H^*\Omega \SU(2) \cong \Z[\gamma^*c]$ similarly. \bigskip

{\bf \S 1.3} The cobordism ring (of suitably prequantized manifolds) $\B_*$ is naturally isomorphic to $\MU_*B\T$. The Hurewicz-Boardman-Chern-Weil Thom characteristic number homomorphism
\[
\hu : \B_*  \cong \MU_*B\T  \to  H_*(\CP_\infty; \MU_*) \to H_*(MU,\Q)[b]   
\]
is injective, and becomes an isomorphism after rationalization. From now on we will identify $\MU_*$ with its image in the symmetric functions, via the Thom isomorphism $H_*\MU \cong H_*B{\rm U} \cong \Sf_*$; this is most naturally expressed \cite{2} in terms of Chern roots with respect to the basis $m$.

Expressions for symplectic cobordism classes such as $[\CP_n ,\omega] := \CP_n(\omega)$ as polynomials in $\omega$ are the topic of \cite{8}; for example $b$ can be naturally identified with $\CP_1(\omega)$. It will be convenient to write $\beta = \CP_1(\omega)/\CP_1$ so $b = \beta \CP_1$, with $\beta$ analogous to $1/kT$. \bigskip

By Quillen, Lazard, and Mi\v{s}\v{c}enko, the dual formal Hopf algebra $\MU^*\CP_\infty \cong \MU^*[[c]]$ is characterized by its group law
\[
c_0 +_\MU c_1 = \exp_\MU (\log_\MU(c_0) + \log_\MU(c_1)) \in \MU^*[[c_0,c_1]] 
\]
with
\[
\log_\MU (c) = \sum \CP_{n-1} \frac{c^n}{n} \;.
\] 
 Let us then write 
\[
\st_\MU(t) = \exp_\G(b \log_\MU(t)) \sim bt + \dots \in (\Sf \otimes \Q[b])[[t]]
\]
by identifying $b\CP_{n-1}$ with the $n$th power sum $p_n$. Ravenel and Wilson \cite{10} show that the generating function 
\[ 
b(t) = 1 + \sum b_{(n)} t^n
\]
is a universal Cartier character \cite{4}
\[
b(t_0) +_\MU b(t_1)  = b(t_0) \cdot b(t_1) \;,
\]
for the Lazard-Quillen ring. \bigskip

${\bf II}$ The Boltzmann exponential

{\bf \S 2.1} A finite subset $E = \{0 < \lambda_1 < \dots < \lambda_l \}$ of the real line defines a step function 
$\lambda : \N \to \R \cup \infty$; these are the energetic sets as in \cite{1}, which define essential examples of canonical ensembles \cite{8} in statistical mechanics. Eigenvalues of Laplacians or spectra of starlight would be examples but that they're infinite \dots

Let 
\[
B : E \mapsto \exp_{BE}(t) := \sum \exp_\G (- \lambda_i t) = - \sum \lambda_i^n \frac{t^n}{n!} 
\]
\[
= - \sum \wp_n(\lambda) \gamma^n(p_1 t)  = - p_1 t - \dots \in (\Sf_* \otimes \Q)[[t]]
\]
define the exponential of the Boltzmann formal group law $BE \in \Q[[X,Y]]$  of the energetic set $E$. Here as usual in the formal world 
\[
\log_\G(x) = - \log |1 - x|, \; \exp_\G = 1 - e^{-x} \in \Q[[x]] \; ,
\]
and the sum on the right is the element of the completed graded algebra of symmetric functions \cite{3} defined by evaluating Newton's power sums 
\[
p_n \mapsto p_n(\lambda) = \sum \lambda^n_i \in \Sf_{2n} \otimes \Q
\]
on (the finite part of) the function $\lambda(i) = \lambda_i.$ The elements $\wp_n = p_1^{-n} p_n$ are their homogenizations or normalizations, eg  $\wp_1 = 1$. Here $p_1$ is a kind of mean value, expected to be invertible; somehow $p_1 \sim 1/kT$ in thermodynamics. 

In this language, the normalized Gibbs free energy 
\[
\Omega(p_1x)  \;  \sim \; - \log_\G \exp_{BE}(x) \in \Q[\wp_*][[p_1x]]
\]
is a generalized cumulant \cite{9}(\S 3.0),\cite{14} generating function, recognizable as a twisted inverse of the fundamental Cartier character of the formal group defined by an energetic ensemble $E$ in the composition
\[
[{\rm Energetic}] \to {\rm Spec} \Sf \times_\Q {\mathbb A}^1 \to \; {\rm Spec} \; \B 
\]
of monoidal functors. 

Cumulants have good additive behavior; I suspect physicists think of them somehow in terms of merging Chomsky clouds of harmonic oscillators.\bigskip

{\bf \S 2.2} Baez's construction extends to define an average for finite ensembles of one-dimensional formal group laws over $\Q$. The symplectic characteristic number homomorphism lifts the classes $\wp_n b^n = p_n$ to define a Gibbs class $\Omega \sim 1 + \sum \Omega_n b^{n+1}  \in \B_* \otimes \Q$, perhap a universal example \cite{6} for the leading (topological) terms in asymptotic expressions such as Stirling's. These constructions have generalizations in noncommutative probability theory \cite{2},\cite{5}(\S 3.3). 

The isomorphism 
\[
H^*_\Q|[X//G]| =  H^0(G,H^*_\Q(X))
\]
(for a finite group $G$ acting on a space $X$) suggests identifying ${\rm MU}_*\otimes \Q$  with the algebra of (naive) characteristic numbers of orbifolds, regarded as the free commutative algebra over $\Q$ on the divided powers
\[
\SP^n \CP(1) := \CP(1)^n/\Sigma_n \cong \CP(n)
\]
of the two-sphere, yielding the Todd class\dots\bigskip

{\bf Acknowledgements} Thanks to K Vogtmann {\tt arXiv:2202.08739} and J Baez for very helpful corespodnence. \bigskip

\bibliographystyle{amsplain}

\end{document}